\documentclass[draft]{amsart}
\usepackage{amssymb}
\usepackage{url}
\usepackage{amscd}
\usepackage{graphicx}
\input epsf
 
 \theoremstyle{plain}
 \newtheorem{theorem}{Theorem}

 \newtheorem{lemma}[theorem]{Lemma}
 
 \theoremstyle{remark}
 \newtheorem{remark}{Remark}
 
 \numberwithin {equation}{section}
 
\begin{document}
%\date{Received ?; received in final form ?}
\title[ Extension Theorems for paraboloids]{ Extension Theorems for Paraboloids in the Finite Field Setting }
\author{Alex Iosevich and Doowon Koh}

\address{Mathematics Department\\
202 Mathematical Sciences Bldg\\
University of Missouri\\
Columbia, MO 65211 USA}
\email{iosevich@math.missouri.edu}
\address{Mathematics Department\\
202 Mathematical Sciences Bldg\\
University of Missouri\\
Columbia, MO 65211 USA}

\email{koh@math.missouri.edu}

%\thanks{Partially supported by grant }

%\subjclass{}

\begin{abstract}
In this paper we study the $L^p-L^r$ boundedness of the extension operators associated with paraboloids in vector spaces over finite fields.
In higher even dimensions, we estimate the number of additive quadruples in the subset $E$ of the paraboloids, that is the number of
quadruples $(x,y,z,w) \in E^4$ with $x+y=z+w.$ As a result, in higher even dimensions, we improve upon the standard
Tomas-Stein exponents (See (\ref{TSexponents}) below) which
Mockenhaupt and Tao (\cite{MT04}) obtained for the boundedness of extension operators for paraboloids 
by estimating the decay of the Fourier transform of measures on paraboloids. In particular, Theorem \ref{p4} below gives the sharp $L^p-L^4$ bound 
up to endpoints in higher even dimensions. Moreover, using the Theorem \ref{p4} and the Tomas-Stein argument, we also study the $L^2-L^r$ estimates.
In the case when $-1$ is not a square number in the underlying finite field, we also study the $L^p-L^r$ bound in higher odd dimensions. 
The discrete Fourier analytic machinery and Gauss sum estimates make an important role in the proof.  

\end{abstract}

\maketitle

%%%%%%%%%%%%%%%%%%%%%%%%%%%%%%%%%%%%%%%%%%%%%%%%%%%%%%%%%%%%%%%%%%%%%%%%%
% Macros
%%%%%%%%%%%%%%%%%%%%%%%%%%%%%%%%%%%%%%%%%%%%%%%%%%%%%%%%%%%%%%%%%%%%%%%%%
 
\newcommand\sfrac[2]{{#1/#2}}
\newcommand\cont{\operatorname{cont}}
\newcommand\diff{\operatorname{diff}}
\tableofcontents
%%%%%%%%%%%%%%%%%%%%%%%%%%%%%%%%%%%%%%%%%%%%%%%%%%%%%%%%%%%%%%%%%%%%%%%%%
 
\section{Introduction}
Let $S$ be a hypersurface in ${\mathbb R}^d$ and $d\sigma$ a surface measure on $S$.
In Euclidean space, the classical extension problems ask one to find the set of exponents $p$ and $r$
such that the estimate 
$$\| (fd\sigma)^\vee \|_{L^r({\mathbb R}^d)} \leq C \|f\|_{L^p(S,d\sigma)} \quad \mbox{for all} \quad f\in L^p(S,d\sigma) $$
holds, where the constant $C>0$ depends only on the exponents $p$ and $r$, and  $(fd\sigma)^\vee $ denotes the inverse Fourier transform of the measure
$fd\sigma.$ The extension problems have received much attention in the last few decades, because they are related to many interesting problems in harmonic 
analysis such as kakeya problems. For a survey of the development of ideas and some recent results on the extension problem, see
\cite{Ta03}. See also \cite{Fe70}, \cite{Zy74}, \cite{Str77}, \cite{Bo91},  \cite{Ste93},   \cite{CI98},  \cite{Ta99}, and \cite{Wo03}, 
and the references therein on recent progress related to this problem. In recent years, the extension problems in the finite field setting have been studied, 
because the finite field case serves as a good model for the Euclidean case. See, for example, \cite{MT04}, \cite{IK07}, and \cite{IKo07}.
Mockenhaupt and Tao (\cite{MT04}) first posed the extension problems in the finite field setting for various algebraic varieties $S$ and 
they obtained reasonably good results from extension problems for paraboloids in $d$-dimensional vector spaces over finite fields.
In two dimensions, the extension problems for the parabola were completely solved by them. Moreover, they  obtained the standard
Tomas-Stein exponents for the $L^p-L^r$ boundedness of extension operators related to paraboloids in higher $d$-dimensional vector spaces over finite fields.
In particular, they improved on the standard Tomas-Stein exponents for the paraboloids in three dimensional vector spaces over finite fields 
in the case when $-1$ is not a square number in the underlying finite fields. The aforementioned authors used the combinatorial
methods to prove the incidence theorems between lines and points in two dimensial vector spaces over finite fields.  
As a result, they could improve upon the standard Tomas-Stein exponents for 
paraboloids in three dimensions. Here and throughout this paper, the standard Tomas-Stein exponents are 
defined as the exponents $p\ge 1$ and $r\ge 1$ such that 
\begin{equation}\label{TSexponents} r \ge \frac{2d+2}{d-1} \quad \mbox{and} \quad r\ge \frac{p(d+1)}{(p-1)(d-1)},\end{equation}
where $d$ denotes the dimension of the vector space over the finite field. 
However, if $-1$ is a square number or the dimension $d$ is greater than three, the combinatorial methods may not be useful
tool to prove the incidence theorems. The purpose of this paper is to significantly improve  the Tomas-Stein exponents for the boundedness of extension operators
associated with paraboloids in higher even dimensional vector spaces over finite fields. Moreover, our results in even dimensions hold without the assumption
that $-1$ is not a square number in the underlying finite field. In higher odd dimensions, we also investigate the cases when the standard 
Tomas-Stein exponents can be improved.  

In order to state our results in detail, let us review some notation and Fourier analytic machinery in the finite field setting.
We denote by ${\mathbb F_q}$ a finite field with $q$ elements whose characteristic is greater than two. Let ${\mathbb F}_q^d$ be a $d$-dimensional 
vector space over the finite field $F_q$. Let $({\mathbb F}_q^d , dx)$ denote a $d$-dimensional vector space, which we endow
with the normalized counting measure $dx$, and $({\mathbb F}_q^d , dm)$ denotes the dual space, which we endow with the counting measure $dm$.
For any complex-valued function $f$ on $({\mathbb F}_q^d , dx), d\ge 1,$ we define the Fourier transform of $f$ by the formula
$$\widehat{f}(m)= q^{-d} \sum_{x\in {\mathbb F}_q^d}\chi(-x\cdot m) f(x)$$
where $\chi$ denotes a non-trivial additive character of ${\mathbb F_q}$, and $ x\cdot m$ is the usual dot product of $x$ and $m$.
Similarly, we define the inverse Fourier transform of the measure $fd\sigma$ by the relation
$$ (fd\sigma)^\vee (m) = \frac{1}{|S|} \sum_{x\in S}\chi(x\cdot m) f(x)$$
where $|S|$ denotes the number of elements in an algebraic variety $S$ in $({\mathbb F}_q^d , dx)$, and $d\sigma$ denotes the normalized
surface measure on $S$. In other words, $ \sigma(x) = \frac{q^d}{|S|} S(x).$
Here and throughout the paper, we identify sets with their characteristic functions. For example, we write $E(x)$ for $\chi_E (x).$
For $1\le p,r <\infty$, we also define
$$\|f\|_{L^p \left({\Bbb F}_q^d,dx \right)}^p=q^{-d}\sum_{x\in {\Bbb 
F}_q^d}|f(x)|^p,$$
$$\|\widehat{f}\|_{L^r \left({\Bbb F}_q^d,dm \right)}^r = \sum_{m\in
{\Bbb F}_q^d}|\widehat{f}(m)|^r$$
 and
$$ \|f\|_{L^p \left(S,d\sigma \right)}^p=\frac{1}{|S|}\sum_{x\in S}|f(x)|^p.$$
Similarly, denote by  $\|f\|_{L^\infty}$ the maximum value of $f$.

\section{ The extension problem and statement of results}
Let $S$ be an algebraic variety in $({\mathbb F}_q^d , dx)$ and $d\sigma$ the surface measure on $S$.
For $1\leq p,r \leq \infty$, we define $ R^*(p\to r)$ by the smallest constant such that the extension estimate 
\begin{equation}\label{DefR^*} \|(fd\sigma)^\vee \|_{L^r({\mathbb F}_q^d , dm)} \leq R^*(p\to r) \|f\|_{L^p(S, d\sigma)}\end{equation}
holds for all functions $f$ on $S$. By duality, $R^*(p\to r)$ is also given by the smallest constant such that the restriction estimate
$$ \|\widehat{g}\|_{L^{p^\prime}(S, d\sigma)} \leq R^*(p\to r) \|g\|_{L^{r^\prime}({\mathbb F}_q^d , dm)} $$
holds for all functions $g$ on $ ({\mathbb F}_q^d , dm).$ In the finite field setting, the extension problem for $S$ asks us to determine the 
set of exponents $p$ and $r$ such that $ R^*(p\to r)\leq C_{p,r} ,$ where the constant $ C_{p,r}$ is independent of $q$ ,
the size of the underlying finite field $F_q$. Using  H\"older's inequality and the nesting properties of $L^p$-norms, we see that
\begin{equation}\label{smallgood}R^*(p_1\to r) \le R^*(p_2\to r) \quad \mbox{for}\quad \ p_1\ge p_2 \end{equation}
and
\begin{equation}\label{rsmallgood}R^*(p\to r_1)\le R^*(p\to r_2) \quad \mbox{for}\quad r_1\ge r_2 \end{equation}
which will allow us to reduce the analysis below to certain endpoint estimates.
We define the paraboloid $S$ in $({\mathbb F}_q^d , dx)$ by the set 
\begin{equation}\label{Defparaboloid} 
S = \{ (\underline{x} , x_d )\in {\mathbb F}_q^d : \underline{x} \in {\mathbb F}_q^{d-1} , 
x_d=\underline{x} \cdot \underline{x} \in {\mathbb F}_q\} .\end{equation}
Our results shall be stated in the following subsections.

\subsection{ $L^p-L^4$ estimate}
\begin{theorem}\label{p4} Let $S$ be the paraboloid in ${\mathbb F}_q^d$ defined as in (\ref{Defparaboloid}).
If $d\geq 4$ is even, then we have 
$$R^*(p\to 4) \lessapprox 1 \quad \mbox{for all} \quad p\geq \frac{4d}{3d-2} $$
\end{theorem}
\begin{remark} Recall that $X \lesssim Y$ means that there exists $C>0$, independent of
$q$ such that $X \leq CY$, and $X \lessapprox Y$, with the controlling parameter $q$, means that for every $\varepsilon >0$ 
there exists $C_{\varepsilon} >0$ such that $X \leq C_{\varepsilon} q^\varepsilon Y.$  
\end{remark}
In higher even dimensions, we claim that Theorem \ref{p4} improves  the standard Tomas-Stein exponents which Mockenhaupt and Tao obtained in \cite{MT04}.
To see this, note that  the standard Tomas-Stein exponents in (\ref{TSexponents}) imply that $ R^*(\frac{4d-4}{3d-5}\to 4) \lesssim 1.$
Since $\frac{4d}{3d-2}$ is less than $\frac{4d-4}{3d-5}$, the claim immediately follows from the inequality in (\ref{smallgood}).
Moreover, Theorem \ref{p4} is sharp  except for the endpoint in the sense that for every $\varepsilon >0, R^*(\frac{4d}{3d-2}-\varepsilon \to 4) \lesssim 1$ fails
to be true.
In fact, the sharpness of Theorem \ref{p4} follows from the following necessary conditions for $R^*(p\to r)\lesssim 1$ in even dimensions $d\geq 4$:
\begin{equation}\label{necessary} r\ge \frac{2d}{d-1}\quad \mbox{and} \quad r \ge \frac{p(d+2)}{(p-1)d}.\end{equation}
In order to prove the necessary conditions, first observe that 
$$ \sum_{m\in {\mathbb F}_q^d} \left|(fd\sigma)^\vee(m) \right|^2 =\frac{q^d}{|S|} \|f\|^2_{L^2(S, d\sigma)}.$$
Using this and Cauchy-Schwartz inequality, we see that if $r\geq 2$ then
$$\frac{q^d}{|S|} \|f\|^2_{L^2(S, d\sigma)} \leq \left(\sum_{m\in {\mathbb F}_q^d} 
\left|(fd\sigma)^\vee(m)\right|^r \right)^\frac{2}{r} \cdot q^{\frac{d(r-2)}{r}} .$$
It therefore follows that
$$ \frac{q^d}{|S|} q^{\frac{-d(r-2)}{r}} \|f\|^2_{L^2(S,d\sigma)} \leq \| (fd\sigma)^\vee \|^2_{L^r({\mathbb F}_q^d, dm)}.$$
Choosing $f$ such that $ \|f\|_{L^2(S,d\sigma)}=1$ and using the fact that $|S|=q^{d-1},$
 we see that we must have
$$ q^{\frac{-d(r-2)}{r}}q \lesssim 1,$$
from which the first estimate in (\ref{necessary}) follows.
On the other hand, if the paraboloid $S$ contains a subspace $H$ of dimension $n$ $(|H|=q^n),$ then we can test 
(\ref{DefR^*}) with $f$ equal to the characteristic function on the set $H$. As a result, we must have that
\begin{equation}\label{necessary2} |H||S|^{-1} q^{\frac{d-n}{r}} \lesssim ( |H||S|^{-1})^\frac{1}{p}.\end{equation}
Since $|H|=q^n$ and $|S|=q^{d-1}$, the inequality in (\ref{necessary2}) implies that 
\begin{equation}\label{generalone}
r\geq \frac{p(d-n)}{(p-1)(d-n-1)}.
\end{equation}
However, if $-1 \in {\mathbb F}_q$ is a square number, say $i^2=-1$ for some $i \in {\mathbb F}_q$, then 
the paraboloid $S$ contains the subspace $H \subset {\mathbb F}_q^d $ of dimension $\frac{d-2}{2}$ which is given by
$$ \left\{ ( s_1, is_1, \ldots, s_k, is_k, \ldots, s_{\frac{d-2}{2}}, is_{\frac{d-2}{2}}, 0, 0 )  :
s_k \in {\mathbb F}_q,\, k=1,2,\ldots, \frac{d-2}{2} \right\}.$$ Since $|H|=q^{\frac{d-2}{2}}$, replacing $n$ in (\ref{generalone})
by $\frac{d-2}{2}$, we obtain the second estimate in (\ref{necessary}).

\begin{figure}[h]
\centering\leavevmode\epsfysize=4.5cm \epsfbox{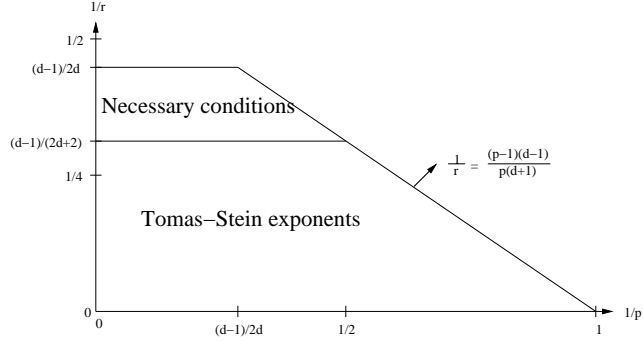}
\caption{\label{Figure1}  In odd dimensions $d\geq 3$, the necessary conditions for $R^*(p\to r)$ bound and the Tomas-Stein exponents.}
\end{figure}

\begin{remark} \label{Re2}
In the case when $d=3$ and $-1$ is not a square number, Mockenhaupt and Tao (\cite{MT04}) obtained an improvement of 
the $"p"$ index of the Tomas-Stein exponent $R^*(2\to 4)$ by showing that $ R^*(\frac{8}{5} \to 4) \lessapprox 1.$
However, if $-1$ is allowed to be a square number ,  one can not expect the improvement of the $"p"$ index of 
the Tomas-Stein exponent $R^*(p\to 4)$ in the case when the dimension $d\geq 3$ is odd. This follows from the fact
that if $-1$ is a square number and $d$ is odd then the paraboloid $S$ always contains the subspace $H \in {\mathbb F}_q^d$ of 
dimension $\frac{d-1}{2}$, defined by 
$$H= \left\{ ( s_1, is_1, \ldots, s_k, is_k, \ldots, s_{\frac{d-1}{2}}, is_{\frac{d-1}{2}}, 0 )  :
s_k \in {\mathbb F}_q,\, k=1,2,\ldots, \frac{d-1}{2} \right\}.$$ 
Note that $|H|= q^{\frac{d-1}{2}}.$ Thus
substituting $n$ in (\ref{generalone}) for $\frac{d-1}{2}$,
a necessary condition for $R^*(p\to r) \lesssim 1$ takes the form
$$  \quad r \ge \frac{p(d+1)}{(p-1)(d-1)}.$$
Note that if $r=4$ then this necessary condition exactly matches the Tomas-Stein exponent in (\ref{TSexponents}) (See Figure \ref{Figure1}).
\end{remark} 
 
\subsection{$L^2-L^r$ estimate}
\begin{figure}[h]
\centering\leavevmode\epsfysize=4.5cm \epsfbox{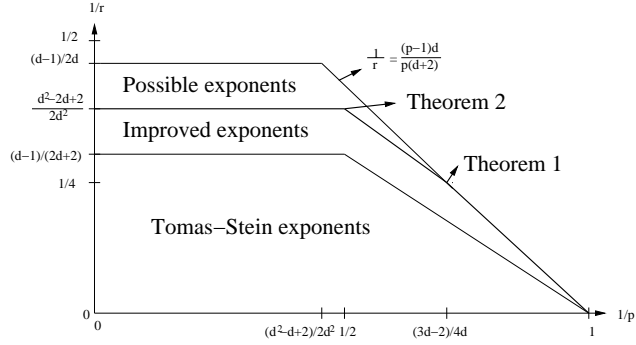}
\caption{\label{Figure2} In even dimensions $d\geq 4,$  exponents for  $R^*(p\to r)$ bound.}
\end{figure}

\begin{theorem}\label{2r}
Let $S$ be the paraboloid in ${\mathbb F}_q^d$ defined as in (\ref{Defparaboloid}).
If $d\geq 4$ is even, then we have 
$$R^*(2\to r) \lessapprox 1 \quad \mbox{for all} \quad r\geq \frac{2d^2}{d^2-2d+2} $$
\end{theorem}
Note that the standard Tomas-Stein exponents in (\ref{TSexponents}) yields that $R^*(2\to \frac{2d+2}{d-1}) \lesssim 1.$
In  higher even dimensions, we therefore see that Theorem \ref{2r} gives an improvement of the "r" index of 
the standard Tomas-Stein exponent $R^*(2\to \frac{2d+2}{d-1}),$ which is true from (\ref{rsmallgood}), 
because the exponent $\frac{2d^2}{d^2-2d+2}$ is less than $\frac{2d+2}{d-1}$ (See Figure \ref{Figure2}).

\begin{remark} Interpolating the results of Theorem \ref{p4} and Theorem \ref{2r} and the trivial bound $R^*(1\to \infty)\lesssim 1$
together with the fact in (\ref{smallgood}), our results can be described as in Figure \ref{Figure2}.
\end{remark}

\subsection{ $L^p-L^r$ estimates in odd dimensions}
As mentioned in Remark \ref{Re2}, in general it is impossible to improve the $"p"$ index of the Tomas-Stein exponent $R^*(p\to 4)$ in odd dimensions.
However, if we assume that $-1$ is not a square number in the underlying finite field ${\mathbb F}_q$, then the improvement of the $"p"$ index of
the Tomas-Stein exponent $ R^*(p\to 4)$ can be obtained in odd dimensions. For example, Mockenhaupt and Tao (\cite{MT04}) obtained the improvement in three dimensions
(see Remark \ref{Re2} ). We shall extend their work to higher odd dimensions $d \ge 7$ in the  specific case.

\begin{theorem}\label{oddimprove}
Let $S$ be the paraboloid in ${\mathbb F}_q^d$ defined as before.
Suppose that $d\ge 7$ is odd and $ q = \underline{p}^l$ for some odd prime $ \underline{p} \equiv 3 \,\,(mod 4) .$
If $l(d-1)$ is not a multiple of four then we have
\begin{equation}\label{Ju1}
R^*(p\to 4) \lessapprox 1 \quad \mbox{for all} \quad p \geq \frac{4d}{3d-2} \end{equation}
and 
\begin{equation}\label{Ju2}
R^*(2\to r) \lessapprox 1 \quad \mbox{for all} \quad r \geq \frac{2d^2}{d^2-2d+2}. \end{equation} 
\end{theorem}

\begin{remark} Observe that the assumptions in Theorem \ref{oddimprove} imply that if   
$l$ is an odd number, $q = \underline{p}^l$ for some odd prime $ \underline{p} \equiv 3 \,\,(mod 4),$ 
and $d=4k+3$ for some $k\in\mathbb {N}$ 
then the conclusions in Theorem \ref{oddimprove} hold. However, since any odd power of the prime $ \underline{p}\equiv 3 \,\,(mod 4)$
is $\equiv 3 \,\,(mod 4),$  we need the condition that $q\equiv 3\,\,(mod 4)$ which means
$-1$ is not a square number in ${\mathbb F}_q.$ In conclusion, we can say that if 
$-1$ is not a square number in ${\mathbb F}_q$ and $d=4k+3$ for some $k\in\mathbb {N}$, then 
the conclusions in Theorem \ref{oddimprove} hold.
 \end{remark}

\subsection{Outline of this paper} In section 3, we shall introduce simple facts related to Gauss sums which 
will be used to prove our main results.
 The proofs of Theorem \ref{p4} and Theorem \ref{2r}  will be given in section 4 and section 5 respectively. 
In the last section, we shall prove Theorem \ref{oddimprove}.
\section{The classical exponential sums}
In this section, we shall collect the well-known facts which follow from estimates of Gauss sums.
Let $\chi$ be a non-trivial additive character of ${\mathbb F}_q$ and $\eta$ a multiplicative character of 
${\mathbb F}_q$ of order two, that is, $ \eta(ab)=\eta(a)\eta(b)$ and  $\eta^2(a)=1$ for all $a,b \in {\mathbb F}_q^*$ but $\eta\not\equiv 1.$
For each $a\in {\mathbb{F}_q}$, the Gauss sum $G_a(\eta, \chi)$ is defined by 
$$ G_a(\eta, \chi) = \sum_{s\in {\mathbb F}_q^*} \eta(s) \chi(as).$$ The magnitude of the Gauss sum is given by the relation
$$ |G_a(\eta, \chi)| = \left\{\begin{array}{ll} q^{\frac{1}{2}} \quad &\mbox{if} \quad a\ne 0\\
                                                  0 \quad & \mbox{if} \quad a=0. \end{array}\right.$$
\begin{remark}
Here, and throughout this paper, we denote by $\chi$ and $\eta$ the canonical additive character  and the quadratic character of 
${\mathbb F}_q$ respectively.
\end{remark}
The following theorem tells us the explicit value of the Gauss sum $G_1(\eta, \chi).$ 
For the nice proof, see \cite{LN97}.
\begin{theorem}\label{ExplicitGauss}
Let ${\mathbb F}_q$ be a finite field with $ q= \underline{p}^l$, where $\underline{p}$ is an odd prime and $l \in {\mathbb N}.$
Then we have
$$G_1(\eta, \chi)= \left\{\begin{array}{ll}  {(-1)}^{l-1} q^{\frac{1}{2}} \quad &\mbox{if} \quad \underline{p} =1 \,\,( mod 4) \\
                    {(-1)}^{l-1} i^l q^{\frac{1}{2}} \quad &\mbox{if} \quad \underline{p} =3 \,\,( mod 4).\end{array}\right. $$
 \end{theorem}

In particular, we have
\begin{equation}\label{square} \sum_{s \in {\mathbb F}_q} \chi (a s^2) = \eta(a) G_1(\eta, \chi) \quad \mbox{for any} \quad a \ne 0,\end{equation}
because $\eta$ is the multiplicative character of ${\mathbb F}_q^* $ of order two. For the nice proof for this equality and
the magnitude of Gauss sums, see \cite{LN97} or \cite{IK04}. As the direct application of the equality in (\ref{square}), we have the following estimate.
\begin{lemma}\label{complete}
For $\beta \in {\mathbb F}_q^k$ and $t\ne 0$, we have
$$ \sum_{\alpha \in {\mathbb F}_q^k} \chi( t \alpha \cdot \alpha + \beta \cdot \alpha ) 
= \chi\left( \frac{\|\beta\|_2}{-4t}\right) \eta^k(t) \left(G_1(\eta, \chi)\right)^k,$$
where, here and throughout the paper, $\|\beta\|_2 = \beta \cdot \beta.$
\end{lemma}
\begin{proof}
It follows that
$$\sum_{\alpha \in {\mathbb F}_q^k} \chi( t \alpha \cdot \alpha + \beta \cdot \alpha ) 
              =\prod_{j=1}^{k} \sum_{\alpha_j\in {\mathbb F}_q} \chi( t\alpha_j^2 + \beta_j\alpha_j).$$
Completing the square in $\alpha_j$-variables, changing of variables, $ \alpha_j+\frac{\beta_j}{2t} \to \alpha_j$,
and using the inequality in (\ref{square}), the proof immediately follows.

\end{proof}
We shall introduce the explicit formula of $ (d\sigma)^\vee$, the inverse Fourier transform of the surface measure
related to the paraboloid. 
\begin{lemma}\label{explicit}
Let $S \subset {\mathbb F}_q^d$ be the paraboloid and $d\sigma$ the surface measure on $ S$. For each 
$m=(\underline{m}, m_d) \in {\mathbb F}_q^{d-1}\times {\mathbb F}_q$ , we have
$$ (d\sigma)^{\vee}(m)= \left\{\begin{array}{ll} q^{-(d-1)} \chi \left( \frac{\|\underline{m}\|_2 }{-4m_d}\right)
\eta^{d-1}(m_d) &\left(G_1(\eta, \chi)\right)^{d-1}
 \quad \mbox{if} \quad m_d \ne 0\\
0 \quad &\mbox{if} \quad m_d =0 , \underline{m} \ne \underline{0}\\
1 \quad &\mbox{if} \quad m=(0,\ldots,0).\end{array}\right.$$                                          
\end{lemma}
\begin{proof}
For each $ m=(\underline{m}, m_d) \in {\mathbb F}_q^d$ , we have

\begin{align*} (d\sigma)^\vee (m) &= \frac{1}{|S|} \sum_{x\in S} \chi(m\cdot x)\\
&=q^{-(d-1)} \sum_{\underline{x} \in {\mathbb F}_q^{d-1}} \chi\left( m_d \,\underline{x}\cdot \underline{x} +\underline{m}\cdot \underline{x} \right). 
\end{align*}
If $m=(0,\ldots,0), $ then it is obvious that $(d\sigma)^{\vee}(m) =1.$
If $ m_d=0, \underline{m} \ne \underline{0}$,  the orthogonality relation of the nontrivial character yields that 
$(d\sigma)^{\vee}(m) =0.$ On the other hand, if $m_d \ne 0$, the proof follows from  Lemma \ref{complete}.
\end{proof}

\section{Proof of the $L^p-L^4$ estimate (Theorem \ref{p4})}
Let $S\subset {\mathbb F}_q^d $ be the paraboloid defined as in (\ref{Defparaboloid}).
Using (\ref{smallgood}) and the usual dyadic pigeonholing argument (see \cite{Gr03}),  it is enough to show that
\begin{equation}\label{Extension} 
\|(Ed\sigma)^\vee\|_{L^4 ({\mathbb F}_q^d,dm)} \lesssim \|E\|_{L^{p_0} (S,d\sigma)} , \quad \mbox{for all} \quad E \subset S,\end{equation}
where $p_0=\frac{4d}{3d-2}.$ Expanding  both sizes in (\ref{Extension}), we see that
$$ \|E\|_{L^{p_0} (S,d\sigma)}= \left(\frac{|E|}{|S|}\right)^{\frac{1}{p_0}}$$ 
and
$$\|(Ed\sigma)^\vee\|_{L^4 ({\mathbb F}_q^d,dm)} = \frac{q^{\frac{d}{4}}}{|S|} \left(\Lambda_4(E)\right)^{\frac{1}{4}}$$
where $\Lambda_4(E)= \sum\limits_{\substack{x,y,z,w\in E\\:x+y=z+w}} 1$ . Since $|S|=q^{d-1}$, it therefore suffices to
show that 
\begin{equation}\label{Easyform}
\Lambda_4(E) \lesssim |E|^{\frac{4}{p_0}} q^{3d-4} q^{\frac{-4d+4}{p_0}}\quad \mbox{for all} \quad E \subset S.
\end{equation} 
We shall need the following estimate.
\begin{lemma}\label{key} Let $S$ be the paraboloid in $({\mathbb F}_q^d , dx)$ defined as before.
In addition, we assume that the dimension of $ {\mathbb F}_q^d, d\ge4,$ is even.
If $E$ is any subset of $S$ then we have
$$\Lambda_4(E) \lesssim \min \{ |E|^3, q^{-1}|E|^3+q^{\frac{d-2}{4}}|E|^{\frac{5}{2}} + q^{\frac{d-2}{2}}|E|^2 \}.$$
\end{lemma}

For a moment, we assume  Lemma \ref{key} which will be proved in the following subsection. We return to the proof of Theorem \ref{p4}.
Note that Lemma \ref{key} implies that if $d\ge 4$ is even and $E$ is any subset of the paraboloid $S$, then 
$$ \Lambda_4(E)\lesssim \left\{\begin{array}{ll}  q^{-1}|E|^3 \quad &\mbox{if} \quad q^{\frac{d+2}{2}} \lesssim |E| \lesssim q^{d-1} \\
                q^{\frac{d-2}{4}}|E|^{\frac{5}{2}} \quad &\mbox{if} \quad q^{\frac{d-2}{2}} \lesssim |E|\lesssim q^{\frac{d+2}{2}} \\
                 |E|^3 \quad &\mbox{if} \quad 1 \lesssim |E| \lesssim q^{\frac{d-2}{2}}.\end{array}\right. $$
Using these upper bounds of $\Lambda_4(E)$ depending on the size of the subset of $S$, the inequality in (\ref{Easyform}) follows by the direct calculation.
Thus the proof of Theorem \ref{p4} is complete.

\subsection{Proof of Lemma \ref{key}}
To prove Lemma \ref{key}, we first note that $ \Lambda_4(E) \le |E|^3$ for all $E \subset S,$ because if we fix
$x,y,z \in E$ then there is at most one $w$  with $x+y=z+w.$ It therefore suffices to show that
\begin{equation}\label{aim1}
\Lambda_4(E) \lesssim q^{-1}|E|^3 + q^{\frac{d-2}{4}}|E|^{\frac{5}{2}} + q^{\frac{d-2}{2}}|E|^2, \quad \mbox{for all} \quad E \subset S.
\end{equation}
Since the set $E$ is a subset of the paraboloid $S$, we have
$$ \Lambda_4(E) =\sum_{\substack{x,y,z,w \in E \\ : x+y=z+w}} 1 \le \sum_{\substack{x,y,z\in E \\ : x+y-z \in S}} 1.$$
For each $x,y,z \in E$, write $ x+y-z = (\underline{x} +\underline{y} -\underline{z}, \,\, x_d+y_d-z_d).$ Then we see that
$x+y-z \in S  \iff \underline{x}\cdot \underline{y} -\underline{y}\cdot \underline{z} -\underline{x}\cdot \underline{z} + \underline{z}\cdot \underline{z} =0,
$ because $x,y,z$ are elements of the paraboloid $S$. It therefore follows that
$$ \Lambda_4(E) \le \sum_{\underline{x},\underline{y},\underline{z} \in \underline{E}} 
\delta_0( \underline{x}\cdot \underline{y} -\underline{y}\cdot \underline{z} -\underline{x}\cdot \underline{z} + \underline{z}\cdot \underline{z})$$
where $\underline{E} =\{ \underline{x}\in {\mathbb F}_q^{d-1} : (\underline{x} , \underline{x}\cdot \underline{x} ) =x \in E \}$, and 
$\delta_0(t) =1$ if $t=0$ and $0$ otherwise. It follows that
\begin{align*}
\Lambda_4(E) &\le \sum_{\underline{x},\underline{y},\underline{z} \in \underline{E}} q^{-1} 
\sum_{s \in {\mathbb F}_q} \chi\left(s( \underline{x}\cdot \underline{y} -\underline{y}\cdot 
\underline{z} -\underline{x}\cdot \underline{z} + \underline{z}\cdot \underline{z})\right) \\
            &= q^{-1} |\underline{E}|^3 + R(\underline{E} )\end{align*}

where $R(\underline{E} )=\sum\limits_{\underline{x},\underline{y},\underline{z} \in \underline{E}} q^{-1} 
\sum\limits_{s \ne 0} \chi\left(s( \underline{x}\cdot \underline{y} -\underline{y}\cdot 
\underline{z} -\underline{x}\cdot \underline{z} + \underline{z}\cdot \underline{z})\right).$
In order to prove the inequality in (\ref{aim1}), it therefore suffices to show that
\begin{equation}\label{aim2}
|{R(\underline{E} )}|^2 \lesssim q^{\frac{d-2}{2}} |\underline{E}|^5 + q^{d-2} |\underline{E} |^4,
\end{equation}
because $|E| =|\underline{E}|.$ Let us estimate ${R(\underline{E} )}^2 .$
The Cauchy-Schwartz inequality applied to the sum in the variable $\underline{x}$ yields
$$ |{R(\underline{E} )}|^2 \le q^{-2} |\underline{E}| \sum_{\underline{x} \in \underline{E} } 
\left| \sum_{\underline{y}, \underline{z}\in \underline{E}, s\ne 0} \chi\left(s( \underline{x}\cdot \underline{y} -\underline{y}\cdot 
\underline{z} -\underline{x}\cdot \underline{z} + \underline{z}\cdot \underline{z})\right) \right|^2.$$
Applying the Cauchy-Schwartz inequality to the sum in the variable $\underline{z}$ and then dominating the sum over $\underline{z} \in \underline{E}$
by the sum over $\underline{z} \in {\mathbb F}_q^{d-1}$, we have

\begin{align*} |{R(\underline{E} )}|^2 \le &q^{-2} |\underline{E}|^2 \sum_{\underline{x} \in \underline{E} } \sum_{\underline{z} \in {\mathbb F}_q^{d-1} }
\left| \sum_{\underline{y}\in \underline{E}, s\ne 0} \chi\left(s( \underline{x}\cdot \underline{y} -\underline{y}\cdot 
\underline{z} -\underline{x}\cdot \underline{z} + \underline{z}\cdot \underline{z})\right) \right|^2\\
=& q^{-2} |\underline{E}|^2 \sum_{\underline{x}\in \underline{E}} M(\underline{x})\end{align*}
where $M(\underline{x})= \sum\limits_{\underline{z} \in {\mathbb F}_q^{d-1} }
\left| \sum\limits_{\underline{y}\in \underline{E}, s\ne 0} \chi\left(s( \underline{x}\cdot \underline{y} -\underline{y}\cdot 
\underline{z} -\underline{x}\cdot \underline{z} + \underline{z}\cdot \underline{z})\right) \right|^2.$
To prove the inequality in (\ref{aim2}), it is enough to show that 
\begin{equation}\label{aim3}
M(\underline{x}) \lesssim q^{\frac{d+2}{2}}|\underline{E}|^2 + q^d |\underline{E}| \quad \mbox{for all} \quad \underline{x} \in \underline{E}.
\end{equation}
Let us estimate the value $M(\underline{x}) $ which is written by 
\begin{align*}
&\sum_{\substack{\underline{z} \in {\mathbb F}_q^{d-1},\\ \underline{y}, \underline{y^\prime} \in \underline{E}, \\ s,s^\prime \ne 0}}
\chi\left(s( \underline{x}\cdot \underline{y} -\underline{y}\cdot 
\underline{z} -\underline{x}\cdot \underline{z} + \underline{z}\cdot \underline{z})\right) 
\chi\left(-s^\prime( \underline{x}\cdot \underline{y^\prime} -\underline{y^\prime}\cdot 
\underline{z} -\underline{x}\cdot \underline{z} + \underline{z}\cdot \underline{z})\right)\\
=&\sum_{\substack{\underline{z}\in {\mathbb F}_q^{d-1} ,\\ \underline{y}, \underline{y^\prime} \in \underline{E}, \\ s,s^\prime \ne 0 
:s=s^\prime}}  \chi\left( (s-s^\prime)\underline{z}\cdot \underline{z} + 
(s(-\underline{y}-\underline{x}) +s^\prime(\underline{y^\prime}+ \underline{x})) \cdot \underline{z}\right)
\chi\left( (s\underline{y} -s^\prime \underline{y^\prime})\cdot \underline{x}\right)\\
&+\sum_{\substack{\underline{z}\in {\mathbb F}_q^{d-1} ,\\ \underline{y}, \underline{y^\prime} \in \underline{E}, \\ s,s^\prime \ne 0 
:s\ne s^\prime}}  \chi\left( (s-s^\prime)\underline{z}\cdot \underline{z} + 
(s(-\underline{y}-\underline{x}) +s^\prime(\underline{y^\prime}+ \underline{x})) \cdot \underline{z}\right)
\chi\left( (s\underline{y} -s^\prime \underline{y^\prime})\cdot \underline{x}\right)\\
=& I +II.
\end{align*}

Since $s=s^\prime \ne 0$ in  the term $I$, we have
\begin{align*}
I =& \sum_{\substack{ \underline{y}, \underline{y^\prime} \in \underline{E},\\ s\ne 0}} \sum_{\underline{z}\in {\mathbb F}_q^{d-1}}
\chi\left(s(-\underline{y}+ \underline{y^\prime})\cdot \underline{z}\right) \chi\left( s(\underline{y} -\underline{y^\prime})\cdot \underline{x} \right)\\
=& q^{d-1} \sum_{\underline{y}\in \underline{E}, s\ne 0} 1 = q^{d-1}(q-1) |\underline{E}| \lesssim q^d |\underline{E}|,
\end{align*}
where the second line follows from the orthogonality relations for the non-trivial character $\chi$ related to the variable $\underline{z}\in {\mathbb F}_q^d.$ 
To complete the proof, it remains to show that 
\begin{equation}\label{aim4} II \lesssim q^{\frac{d+2}{2}} |\underline{E}|^2.\end{equation}
Setting $a=s, b =\frac{s^\prime}{s},$ we see that
$$II = \sum_{\substack{\underline{z} \in {\mathbb F}_q^{d-1},\\ \underline{y}, \underline{y^\prime} \in \underline{E}, \\ a\ne 0 , b \ne 0,1}}
\chi\left(a(1-b)\underline{z}\cdot \underline{z} + a(-\underline{y}-\underline{x}+ b(\underline{y^\prime}+ \underline{x}))\cdot \underline{z} \right)
\chi\left(a(\underline{y}-b\underline{y^\prime})\cdot \underline{x}\right).$$
Control the sum over $\underline{z} \in {\mathbb F}_q^{d-1}$ by using Lemma \ref{complete} and then the term $II$ is given by $(G_1(\eta, \chi))^{d-1}$ times
$$  \sum_{\substack{\underline{y}, \underline{y^\prime} \in \underline{E},\\ a\ne 0, b\ne 0,1}}
\eta(1-b) \eta(a)\chi\left(\left[ \frac{\|(-\underline{y}-\underline{x})+b(\underline{y^\prime}+\underline{x})\|_2}{-4(1-b)} + 
(\underline{y}-b\underline{y^\prime})\cdot \underline{x}\right] a \right) , $$
where we also used the fact that $\eta^{d-1} = \eta$ , because $d$ is even and $\eta$ is a multiplicative character of order two.
Since the sum over the variable $a\in {\mathbb F}_q^*$ is the Gauss sum, we conclude that
$$II \lesssim q^{\frac{d-1}{2}}q^{\frac{1}{2}} (q-2) |\underline{E}|^2 \lesssim q^{\frac{d+2}{2}} |\underline{E}|^2. $$
Thus the inequality in (\ref{aim4}) holds and the proof of Lemma \ref{key} is complete.

\section{Proof of the $L^2-L^r$ estimate (Theorem \ref{2r})}
In order to prove  Theorem \ref{2r}, we shall use Theorem \ref{p4} together with interpolation theorems and some specific properties of
the Fourier transform of the paraboloid which were investigated by Mockenhaupt and Tao (\cite {MT04}). See also \cite{Gr03}.
By (\ref{rsmallgood}), it suffices to show that
$$R^*\left(2\to \frac{2d^2}{d^2-2d+2}\right) \lessapprox 1 \quad \mbox{whenever} \quad d\ge 4\quad \mbox{is even}.$$
We shall use the Tomas-Stein argument.
Let $R^*: L^p(S,d\sigma) \to L^r({\mathbb F}_q^d, dm)$ be the extension map  $f\to (fd\sigma)^\vee$ , and 
$R: L^{r^\prime}({\mathbb F}_q^d, dm)\to L^{p^\prime}(S,d\sigma)$ be its dual, the restriction map  
$g\to \widehat{g}|_S.$
Note that $ R^*Rg = (\widehat{g}d\sigma)^\vee = g \ast (d\sigma)^\vee$ for all function $g$ on ${\mathbb F}_q^d.$
By the Tomas-Stein argument, it therefore suffices to show that 
$$  \| g \ast(d\sigma)^\vee \|_{L^{p_0}({\mathbb F}_q^d, dm)} 
\lessapprox \|g\|_{L^{p_0^\prime}({\mathbb F}_q^d , dm)}$$
where $p_0=\frac{2d^2}{d^2-2d+2}$ and $p_0^\prime = \frac{2d^2}{d^2+2d-2}.$
Note that 
\begin{align*}\| g \ast \delta_0 \|_{L^{p_0}({\mathbb F}_q^d, dm)} &=\| g  \|_{L^{p_0}({\mathbb F}_q^d, dm)} \\
&\leq \| g \|_{L^{p_0^\prime}({\mathbb F}_q^d, dm)} \end{align*}
where the last line follows from the facts that $dm$ is the counting measure and $ p_0 > p_0^\prime.$
Thus it is enough to show that 
\begin{equation}\label{final1}
\| g \ast K\|_{L^{p_0}({\mathbb F}_q^d, dm)} \lessapprox \| g \|_{L^{p_0^\prime}({\mathbb F}_q^d, dm)} ,\end{equation}
where $K$ is the Bochner-Riesz kernel given by $ K = (d\sigma)^\vee -\delta_0.$ 
Recall that $K(m)=0$ if $m=(0,\ldots,0)$, and $ K(m)= (d\sigma)^\vee(m)$ otherwise.
We now claim that the following two estimates hold:
\begin{equation}\label{e1}
\| g \ast K\|_{L^2({\mathbb F}_q^d, dm)} \lessapprox q \| g \|_{L^2({\mathbb F}_q^d, dm)} \end{equation}
and 
\begin{equation}\label{e2}
\| g \ast K\|_{L^4({\mathbb F}_q^d, dm)} \lessapprox q^{\frac{-d+4}{4}} \| g \|_{L^{\frac{4d}{3d-2}}({\mathbb F}_q^d, dm)} .\end{equation}
The inequality in (\ref{e1}) follows from elementary properties of Fourier transform. In fact, we have
\begin{align*}\| g \ast K\|_{L^2({\mathbb F}_q^d, dm)}&= \| \widehat{g} \widehat{K}\|_{L^2({\mathbb F}_q^d, dx)}\\
                                                      &\leq \| \widehat{K}\|_{L^\infty({\mathbb F}_q^d, dx)} \| {g}\|_{L^2({\mathbb F}_q^d, dm)}.\end{align*}
However, $ \widehat{K}(x) =d\sigma(x)- \widehat{\delta_0}(x) = q S(x)-1 \leq q.$ Therefore the inequality in (\ref{e1}) holds.
For a moment, we assume the inequality in (\ref{e2}) which will be proved in the following subsection.
Note that the operator $g\to g\ast K$ is self-adjoint, because $K$ is essentially the inverse Fourier transform of real-valued function $d\sigma.$
Thus the inequality in (\ref{e2}) implies that 
\begin{equation}\label{e3}\| g \ast K\|_{L^{\frac{4d}{d+2}}({\mathbb F}_q^d, dm)} 
\lessapprox q^{\frac{-d+4}{4}} \| g \|_{L^{\frac{4}{3}}({\mathbb F}_q^d, dm)} .\end{equation}
Interpolating (\ref{e2}) and (\ref{e3}) (with $\theta =\frac{1}{2}$ in the Riesz-Thorin theorem), we also see that
\begin{equation}\label{e4}
\| g \ast K\|_{L^{\frac{4d}{d+1}}({\mathbb F}_q^d, dm)} 
\lessapprox q^{\frac{-d+4}{4}} \| g \|_{L^{\frac{4d}{3d-1}}({\mathbb F}_q^d, dm)} .\end{equation}
Once again interpolate $(\ref{e1})$ and $(\ref{e4})$ ( with $\theta=\frac{4}{d}$ in the Riesz-Thorin theorem),
and then we obtain the inequality in (\ref{final1}). Thus the proof is complete.
\subsection{Proof of the estimate in (\ref{e2})}
For $d\ge 4$ even, we must show that
$$ \| g \ast K\|_{L^4({\mathbb F}_q^d, dm)} \lessapprox q^{\frac{-d+4}{4}} \| g \|_{L^{\frac{4d}{3d-2}}({\mathbb F}_q^d, dm)} $$
where $g$ is the arbitrary function on ${\mathbb F}_q^d$ and $K=(d\sigma)^\vee -\delta_0.$ 
Let $m=(m_1, \ldots, m_{d-1}, m_d) \in {\mathbb F}_q^d.$ For each $a\in {\mathbb F}_q$ and the function $g$ on ${\mathbb F}_q^d$,
we define the function $g_a$ as the restriction of $g$ to the hyperplane $ \{m\in {\mathbb F}_q^d : m_d =a\}$.
Then it is enough to show that for each $a\in {\mathbb F}_q,$ the estimate
\begin{equation}\label{s1}
\| g_a \ast K\|_{L^4({\mathbb F}_q^d, dm)} \lessapprox q^{\frac{-d^2+3d-2}{4d}} \| g_a \|_{L^{\frac{4d}{3d-2}}({\mathbb F}_q^d, dm)}\end{equation}
holds, because the estimate $(\ref{s1})$ yields the following estimates:
\begin{align*} \| g \ast K\|_{L^4({\mathbb F}_q^d, dm)} 
&\leq \sum_{a \in {\mathbb F}_q} \| g_a \ast K\|_{L^4({\mathbb F}_q^d, dm)}\\
&\lessapprox q^{\frac{-d^2+3d-2}{4d}}\sum_{a \in {\mathbb F}_q}\| g_a \|_{L^{\frac{4d}{3d-2}}({\mathbb F}_q^d, dm)} \\
&\leq q^{\frac{-d^2+3d-2}{4d}} q^{\frac{d+2}{4d}} 
\left(\sum_{a \in {\mathbb F}_q}\| g_a \|_{L^{\frac{4d}{3d-2}}({\mathbb F}_q^d, dm)}^{\frac{4d}{3d-2}} \right)^\frac{3d-2}{4d}\\
&= q^{\frac{-d+4}{4}} \|g\|_{L^{\frac{4d}{3d-1}}({\mathbb F}_q^d, dm)} . \end{align*}
Without loss of generality, we may assume that $a=0$, because of translation invariance. Thus it suffices to show that
\begin{equation}\label{s2}\| g_0 \ast K\|_{L^4({\mathbb F}_q^d, dm)} 
\lessapprox q^{\frac{-d^2+3d-2}{4d}} \| g_0 \|_{L^{\frac{4d}{3d-2}}({\mathbb F}_q^d, dm)}\end{equation}
By Lemma \ref{explicit} and the definition of $K$, note that 
$$K(\underline{m}, m_d) = q^{-(d-1)}\chi \left( \frac{\|\underline{m}\|_2 }{-4m_d}\right)
\eta^{d-1}(m_d) \left(G_1(\eta, \chi)\right)^{d-1}$$
when $m_d \ne 0$ and $ K(\underline{m},0) = 0.$ Using this and the definition of the function $g_0$, the left-hand side of (\ref{s2}) is given by
\begin{equation}\label{form1}
q^{\frac{-d+1}{2}}\left(\sum_{\underline{m} \in {\mathbb F}_q^{d-1}} \sum_{m_d \ne 0} \left| \sum_{\underline{m^\prime} \in {\mathbb F}_q^{d-1}}
g(\underline{m^\prime}, 0) \chi\left( \frac{\| \underline{m}-\underline{m^\prime} \|_2}{4 m_d} \right) \right|^4 \right)^{\frac{1}{4}},
\end{equation}
where we used the facts that  $|G_1(\eta, \chi)| =q^{\frac{1}{2}}$ and $|\eta|=1$ and then changing of variables, $-m_d \to m_d.$
Letting $\underline{u} = \frac{-\underline{m}}{2m_d}$ and $s= \frac{1}{4m_d}$, we see that
$$ \frac{\| \underline{m}-\underline{m^\prime} \|_2}{4 m_d} 
= \frac{1}{4s} \,\underline{u} \cdot \underline{u} + \underline{u}\cdot \underline{m^\prime} + s\, \underline{m^\prime}\cdot \underline{m^\prime}.$$
Thus the term $(\ref{form1})$ becomes 
\begin{align*}
&q^{\frac{-d+1}{2}}\left(\sum_{\underline{u} \in {\mathbb F}_q^{d-1}} \sum_{s \ne 0} 
\left|\chi(\frac{1}{4s} \,\underline{u} \cdot \underline{u} ) \sum_{\underline{m^\prime} \in {\mathbb F}_q^{d-1}}
g(\underline{m^\prime}, 0) \chi\left(  (\underline{u}, s) \cdot (\underline{m^\prime}, \underline{m^\prime}\cdot \underline{m^\prime})\right) \right|^4 \right)^\frac{1}{4}\\
\leq & q^{\frac{-d+1}{2}} \left( \sum_{(\underline{u}, s) \in {\mathbb F}_q^d} \left|(Gd\sigma)^\vee(\underline{u}, s) \right|^4 \right)^\frac{1}{4} 
= q^{\frac{-d+1}{2}} \|(G d\sigma)^\vee \|_{L^4({\mathbb F}_q^d, dm)}, \end{align*}
where the function  $G$ on the  paraboloid $S$ is defined by 
$$ G(\underline{x}, \underline{x}\cdot \underline{x}) =|S| g(\underline{x}, 0)=q^{d-1}g(\underline{x}, 0).$$
Using Theorem \ref{p4}, this can be bounded by 
$$\lessapprox q^{\frac{-d+1}{2}} \|G\|_{L^{\frac{4d}{3d-2}}(S,d\sigma)}
= q^{\frac{-d^2+3d-2}{4d}} \|g_0\|_{L^{\frac{4d}{3d-2}}({\mathbb F}_q, dm)},$$
where the equality follows from the observation that 
\begin{align*}\|G\|_{L^{\frac{4d}{3d-2}}(S,d\sigma)} =& \left( \frac{1}{|S|} \sum_{x\in S} |G(x)|^\frac{4d}{3d-2}\right)^\frac{3d-2}{4d}\\
=&q^{\frac{-(d-1)(3d-2)}{4d}} \left(\sum_{\underline{m} \in {\mathbb F}_q^{d-1}} \left| q^{d-1} g(\underline{m}, 0)\right|^\frac{4d}{3d-2}\right)^\frac{3d-2}{4d}\\
=& q^{\frac{(d+2)(d-1)}{4d}}\|g_0\|_{L^{\frac{4d}{3d-2}}({\mathbb F}_q, dm)}.\end{align*}
Thus the estimate in (\ref{s2}) holds and so the proof is complete.

\section{Proof of the $L^p-L^r$ estimates in odd dimensions (Theorem \ref{oddimprove})}
Let us prove the first part (\ref{Ju1}) of Theorem \ref{oddimprove}.
Using the arguments for the proof of Theorem \ref{p4}, it suffices to show that
\begin{equation}\label{Doowon1}\Lambda_4(E) \lesssim |E|^{\frac{4}{p_0}} q^{3d-4} q^{\frac{-4d+4}{p_0}}
\quad \mbox{for all} \quad E \subset S , \end{equation}
where $p_0 = \frac{4d}{3d-2}$ and $\Lambda_4(E)= \sum\limits_{\substack{x,y,z,w\in E\\:x+y=z+w}} 1$ .
We shall use the following lemma.

\begin{lemma}\label{samekey} Let $S$ be the paraboloid in ${\mathbb F}_q^d$ defined as before.
Suppose that $d$ is odd and $ q = \underline{p}^l$ for some prime $ \underline{p} = 3 \,\,(mod 4) .$
If $l(d-1)$ is not a multiple of four and $E$ is any subset of $S$, then we have
$$\Lambda_4(E) \lesssim \min \{ |E|^3, q^{-1}|E|^3+q^{\frac{d-3}{4}}|E|^{\frac{5}{2}} + q^{\frac{d-2}{2}}|E|^2 \}.$$ \end{lemma}
We assume this lemma for a moment. By Lemma \ref{samekey}, we see that 
$$ \Lambda_4(E)\lesssim \left\{\begin{array}{ll}  q^{-1}|E|^3 \quad &\mbox{if} \quad q^{\frac{d+1}{2}} \lesssim |E| \lesssim q^{d-1} \\
                q^{\frac{d-3}{4}}|E|^{\frac{5}{2}} \quad &\mbox{if} \quad q^{\frac{d-1}{2}} \lesssim |E|\lesssim q^{\frac{d+1}{2}} \\
                q^{\frac{d-2}{2}}|E|^2 \quad &\mbox{if} \quad q^{\frac{d-2}{2}} \lesssim |E|\lesssim q^{\frac{d-1}{2}} \\
                 |E|^3 \quad &\mbox{if} \quad 1 \lesssim |E| \lesssim q^{\frac{d-2}{2}}.\end{array}\right. $$
Note that this estimate yields the inequality in (\ref{Doowon1}). Thus it is enough to show that Lemma \ref{samekey} holds.
Let us prove Lemma \ref{samekey}. Let $E$ be a subset of the paraboloid $S$. Since $ \Lambda_4(E) \leq  |E|^3$ for each $E\in S$,
it suffice to show 
\begin{equation}
\Lambda_4(E) \lesssim q^{-1}|E|^3+q^{\frac{d-3}{4}}|E|^{\frac{5}{2}} + q^{\frac{d-2}{2}}|E|^2 \quad \mbox{for} \quad E \subset S .
\end{equation}
Define the set $\underline{E} \subset {\mathbb F}_q^{d-1}$ by
$\underline{E} =\{ \underline{x}\in {\mathbb F}_q^{d-1} : (\underline{x} , \underline{x}\cdot \underline{x} ) =x \in E \}.$
 Repeating the same arguments as in the proof of Lemma \ref{key}, we see that
$$\Lambda_4(E) \leq q^{-1} |\underline{E}|^3 + R(\underline{E} ),$$
where $R(\underline{E} )=\sum\limits_{\underline{x},\underline{y},\underline{z} \in \underline{E}} q^{-1} 
\sum\limits_{s \ne 0} \chi\left(s( \underline{x}\cdot \underline{y} -\underline{y}\cdot 
\underline{z} -\underline{x}\cdot \underline{z} + \underline{z}\cdot \underline{z})\right).$ 
From this and the fact that $|E|=|\underline{E}|$, it is enough to show that
$$|{R(\underline{E} )}|^2 \lesssim q^{\frac{d-3}{2}} |\underline{E}|^5 + q^{d-2} |\underline{E} |^4.$$
As in the proof of  Lemma \ref{key}, we have
$$|{R(\underline{E} )}|^2 \leq q^{-2} |\underline{E}|^2 \sum_{\underline{x}\in \underline{E}} M(\underline{x}) ,$$
where $M(\underline{x})= \sum\limits_{\underline{z} \in {\mathbb F}_q^{d-1} }
\left| \sum\limits_{\underline{y}\in \underline{E}, s\ne 0} \chi\left(s( \underline{x}\cdot \underline{y} -\underline{y}\cdot 
\underline{z} -\underline{x}\cdot \underline{z} + \underline{z}\cdot \underline{z})\right) \right|^2.$
It therefore suffices to show that
$$M(\underline{x}) \lesssim q^{\frac{d+1}{2}}|\underline{E}|^2 + q^d |\underline{E}| \quad \mbox{for all} \quad \underline{x} \in \underline{E}.$$
As in the proof of Lemma \ref{key}, we also see that
$$M(\underline{x}) \lesssim q^d |\underline{E}| + II,$$
where $$II = \sum_{\substack{\underline{z} \in {\mathbb F}_q^{d-1},\\ \underline{y}, \underline{y^\prime} \in \underline{E}, \\ a\ne 0 , b \ne 0,1}}
\chi\left(a(1-b)\underline{z}\cdot \underline{z} + a(-\underline{y}-\underline{x}+ b(\underline{y^\prime}+ \underline{x}))\cdot \underline{z} \right)
\chi\left(a(\underline{y}-b\underline{y^\prime})\cdot \underline{x}\right).$$
Thus our final task is to show that
\begin{equation}\label{taget}
II \lesssim q^{\frac{d+1}{2}}|\underline{E}|^2.
\end{equation} 
To estimate $II$, calculate the sum over $\underline{z} \in {\mathbb F}_q^{d-1}$ by using Lemma \ref{complete} and then we see that 
the term $II$ is given by 
$$ (G_1(\eta, \chi))^{d-1} \sum_{\substack{\underline{y}, \underline{y^\prime} \in \underline{E},\\ a\ne 0, b\ne 0,1}}
\chi\left(\left[ \frac{\|(-\underline{y}-\underline{x})+b(\underline{y^\prime}+\underline{x})\|_2}{-4(1-b)} + 
(\underline{y}-b\underline{y^\prime})\cdot \underline{x}\right] a \right) , $$
where we also used the fact that $\eta^{d-1} = 1$ , because $d$ is odd and $\eta$ is a multiplicative character of order two.
Using Theorem \ref{ExplicitGauss} together with the assumptions of Lemma \ref{samekey} , we see that $(G_1(\eta, \chi))^{d-1} =-q^{\frac{d-1}{2}}.$
Letting $\Gamma_{\underline{x}}(\underline{y},\underline{y^\prime}, b) =
\left[ \frac{\|(-\underline{y}-\underline{x})+b(\underline{y^\prime}+\underline{x})\|_2}{-4(1-b)} + 
(\underline{y}-b\underline{y^\prime})\cdot \underline{x}\right]$, it therefore follows that
\begin{align*}II =& -q^{\frac{d-1}{2}} \sum_{\substack{\underline{y}, \underline{y^\prime} \in \underline{E}, \\ a\ne 0 , b \ne 0,1 \\: 
\Gamma_{\underline{x}}(\underline{y},\underline{y^\prime}, b)\ne 0} }
\chi\left( \Gamma_{\underline{x}}(\underline{y},\underline{y^\prime}, b) \cdot a\right) 
-q^{\frac{d-1}{2}} \sum_{\substack{\underline{y}, \underline{y^\prime} \in \underline{E}, \\ a\ne 0 , b \ne 0,1 \\: 
\Gamma_{\underline{x}}(\underline{y},\underline{y^\prime}, b)= 0} }1 \\
\leq & -q^{\frac{d-1}{2}} \sum_{\substack{\underline{y}, \underline{y^\prime} \in \underline{E}, \\ a\ne 0 , b \ne 0,1 \\: 
\Gamma_{\underline{x}}(\underline{y},\underline{y^\prime}, b)\ne 0} }
\chi\left( \Gamma_{\underline{x}}(\underline{y},\underline{y^\prime}, b) \cdot a\right) \\
=&q^{\frac{d-1}{2}} \sum_{\substack{\underline{y}, \underline{y^\prime} \in \underline{E},   b \ne 0,1 \\: 
\Gamma_{\underline{x}}(\underline{y},\underline{y^\prime}, b)\ne 0} } 1 \leq q^{\frac{d+1}{2}} |\underline{E}|^2,
\end{align*}
where the third line follows from the fact that $ \sum\limits_{ a \in {\mathbb F}_q^*} \chi( t\cdot a) =-1 $ for $t\ne 0$.
Thus the inequality in (\ref{taget}) holds and we complete the proof of the first part (\ref{Ju1}) of Theorem \ref{oddimprove}.
The proof of the second part (\ref{Ju2}) of Theorem \ref{oddimprove} immediately follows from the same arguments as in the proof of 
Theorem \ref{2r}.

  \end{document}